\pdfoutput=1
\documentclass[smallcondensed]{svjour3}     %
\usepackage[utf8]{inputenc}
\usepackage{url}
\usepackage[misc]{ifsym}
\usepackage{fixltx2e}
\usepackage{float}
\usepackage{amssymb}
\usepackage[dvips]{graphicx}
\usepackage{color}
\usepackage{subfigure}
\usepackage{amsmath}%

\date{\today}

\newcommand{\setR}{{\mathord{\mathbb R}}}

\newcommand{\Tr}{{\mathrm{Tr}}}

\begin{document}

\title{A nonlinear preconditioner for experimental design
  problems}

\author{Mario S. Mommer \and Andreas Sommer \and Johannes
  P. Schl{\"o}der \and H. Georg Bock}

\institute{Mario S. Mommer (\Letter)\and Andreas Sommer \and Johannes  P. Schl{\"o}der \and H. Georg Bock
\at
Interdisciplinary Center for Scientific Computing (IWR), Heidelberg
University,
INF 368, 69120 Heidelberg, Germany.\\
\email{\tt {\scriptsize \{mario.mommer|andreas.sommer|johannes.schloeder|bock\}@iwr.uni-heidelberg.de}}}

\maketitle

\begin{abstract}
  We address the slow convergence and poor stability of quasi-newton
  sequential quadratic programming (SQP) methods that is observed when
  solving experimental design problems, in particular when they are
  large.  Our findings suggest that this behavior is due to the fact
  that these problems often have bad absolute condition numbers.  To
  shed light onto the structure of the problem close to the solution,
  we formulate a model problem (based on the $A$-criterion), that is
  defined in terms of a given initial design that is to be improved.
  We prove that the absolute condition number of the model problem
  grows without bounds as the quality of the initial design improves.
  Additionally, we devise a preconditioner that ensures that the
  condition number will instead stay uniformly bounded.  Using
  numerical experiments, we study the effect of this reformulation on
  relevant cases of the general problem, and find that it leads to
  significant improvements in stability and convergence behavior.

\end{abstract}

\section{Introduction}

One of the most important aspects in model based optimization and
investigation of real world processes is the estimation of parameters
appearing in a model. Typically, one estimates these parameters
solving a regression problem based on data collected from one or more
experiments. Optimal experimental design is the task of choosing the
best experimental setup from a set of possible ones, and according to
a predefined criterion. As this task is bound to be constrained in
non-trivial ways, and is formulated in terms of the optimality
conditions of an underlying regression problem, it presents a rich
class of challenging optimization problems. It is also a practically
relevant class, as the solution of these problems leads to important
improvements in the efficiency of research and development
\cite{Franceschini20084846,Schoeneberger2008}.

A standard setting is the estimation of parameters using 
weighted least-squares regression. It is
then natural to rate the experiments according to the quality of the
Fisher information matrix, or of its inverse, the variance-covariance
matrix. A variety of criteria to rate this quality exists, and they are
traditionally named after letters of the alphabet. In this article we
will focus mainly on the so called {\em $A$-criterion}, which is
defined as the trace of the variance-covariance matrix of the
regression. For a full list and further discussion we refer to the
classic texts \cite{Pukelsheim2006,Fedorov1972a}.

The experimental design problems we are concerned with here have the
following form. We consider the nonlinear regression problem
\begin{equation}
\label{eq:3}
\tilde{p}=\mathop{\mathrm{arg min}}_{p} \sum_{i=1}^{m} |F(p,q,t_i)-\mu_{i}|^2
\end{equation}
where $F$ is a nonlinear function depending on {\em parameters} $p$,
{\em controls} $q$, and {\em measurement points} $t_i$. Often these
points refer to measurement times, but our scope is not restricted to
this case. The values $\mu_i$ represent the results of measurements,
and $F$ represents the model under study. If the measurement errors are
independent, normally distributed with variance one and zero mean, the
parameters obtained from (\ref{eq:3}) will be a random variable which,
to a first approximation, is drawn from a multivariate normal
distribution centered around the true parameter values $p^{*}$, and
with variance-covariance matrix $\Sigma=(J^{T}J)^{-1}$, where the
matrix $J$ is given by
\[ J_{ij}=\frac{\partial}{\partial p_{j}}F(p^{*},q,t_i). \] The
experiment design problem we consider is the problem of finding
controls $q$, and a subset of measurement points such that the trace
of the variance-covariance matrix is minimal, at least for the
parameters $p$ that we believe to be plausible when designing the
experiment.

Our basic problem is thus, given a set of feasible controls $\Theta$,
and a {\em measurement budget} $m_{\max}<m$, find a set
$M\subset\{1,2,\ldots,m\}$ with cardinality $\# M = m_{max}$, and a
$q\in\Theta$ such that
\[ \mathrm{Tr}\left(\left[J_{[M]}(q)^TJ_{[M]}(q)\right]^{-1}\right) \]
is minimal. Here, we have made explicit the dependence on $q$, and
have denoted by $J_{[M]}$ the matrix $J$ after deleting all rows whose
index is not in the set $M$.

There are a couple of approaches in the literature to solve this
problem
\cite{Koerkel2002fixed,Bauer2000,Lohmann1992a,Fedorov1972a,Franceschini20084846}. One
important approach, and the one we will consider here, consists in
using a {\em relaxed formulation} to obtain a constrained minimization
problem in continuous variables, and then to apply a modern
optimization method, typically in the form of a sequential quadratic
programming (SQP) \cite{Han1977,Powell1978a,Nocedal2006} algorithm to
solve it. This approach was pioneered by
\cite{Lohmann1992a,Bauer2000,Koerkel2002fixed}, and the formulation has the
following form.

Define $W:\setR^{m}\rightarrow\setR^{m\times m}$ through
$W(w):=\mathrm{diag}(w)$, a diagonal matrix whose entries are the
elements of the vector of weights $w$. We define the set of admissible
weights
\[ \Omega(m_{\max},m):=\{w\in[0,1]^m | \sum_{i=1}^mw_i=m_{\max}\}\]
The problem now is
\begin{subequations}
\label{eq:csdproblem}
\begin{align}
&\min_{w,q}\mathrm{Tr}\left(\left[J^{T}(q)W(w)J(q)\right]^{{-1}}\right)
\intertext{subject to} 
&q\in\Theta,\quad w\in\Omega(m_{\max},m).
\end{align}
\end{subequations}
The form of this problem is such that it can be given to an SQP solver
as it is.  It has the important advantage that the optimization of the
controls and of the weights occurs at the same time. On the other
hand, it has the disadvantage of not necessarily yielding integer
solutions, although practical experience shows that either the
solution is integer, or only a few weights are not, which can then be
remedied by using a rounding technique \cite{Pukelsheim1992}. The
understanding of this phenomenon is incomplete, but see
\cite{Sager2012corrected} for an important contribution on this issue.
Alternatively, the weights in (\ref{eq:csdproblem}) can be interpreted
as a required precision of the measurements, i.e. as a measure of the
maximum variance of the error which is acceptable.

\begin{figure}
\includegraphics[width=12cm]{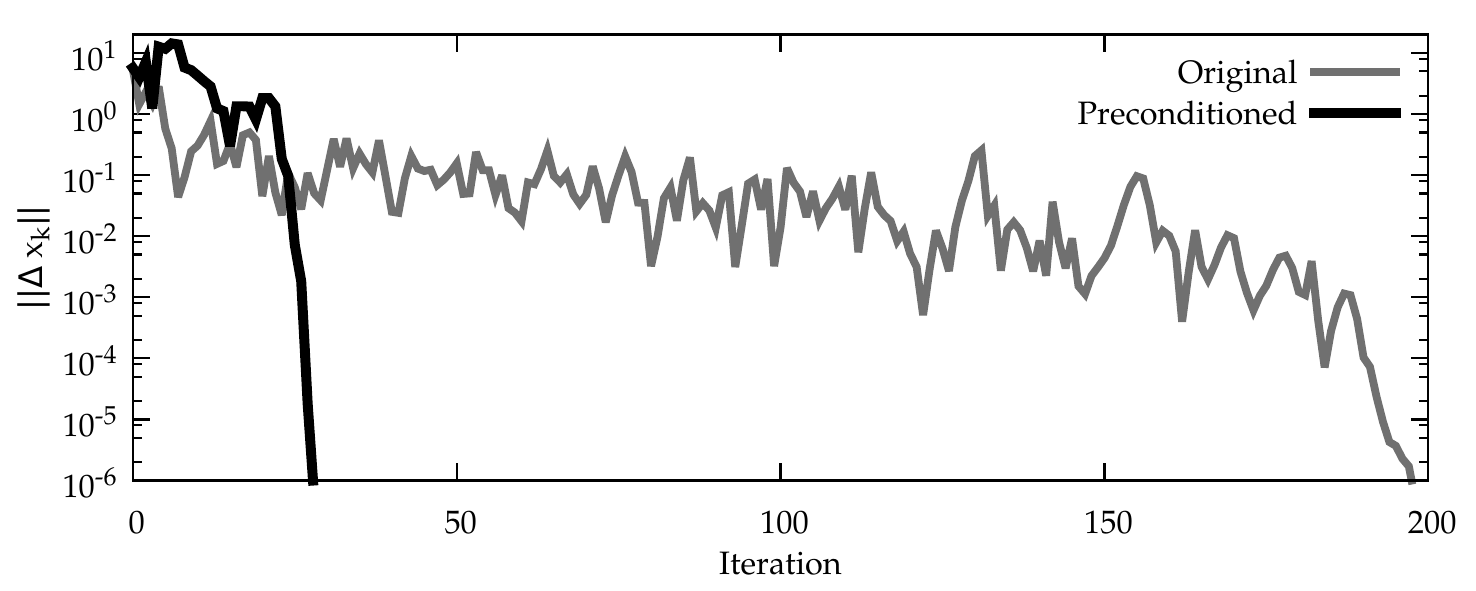}
\caption{Convergence behavior of an SQP algorithm solving a nonlinear
  experimental design problem.}
\label{fig:gure}
\end{figure}

In practice, it turns out that solving the problem
(\ref{eq:csdproblem}) in this way leads to poor convergence. Typical
behavior can be seen in Figure \ref{fig:gure}, where we have plotted
the length of the search direction obtained from the quadratic problem
against the iteration number (the precise description of the numerical
experiment can be found in section \ref{sec:naexps}). It is important
to remark that this difficulty appears also when there are no external
controls, i.e. when the vector of controls $q$ is empty.

In what follows, we will attempt to shed some light on the question as
to why this behavior emerges, and in particular onto how to solve
it. To this end, inspired by the ``test equation''
\cite{Dahlquist1963} used in the theory of stiff ordinary differential
equations, we tackle the generality of the setting by developing a
model problem. In this model problem we discover arbitrarily bad
absolute conditioning under fairly generic conditions. It turns out
that this can be corrected by introducing a simple but nonobvious
transformation. This transformation can also be applied to
(\ref{eq:csdproblem}), yielding a problem for which the SQP method
converges in much less iterations than for the original one.

While the connection between absolute condition numbers and the
difficulty of solving an optimization problem has been made before
\cite{Zolezzi2002,Zolezzi2003}, its role in the convergence behavior
of SQP methods has not, to our knowledge, been investigated.  Thus,
while we cannot provide a complete theoretical justification for the
slow convergence of SQP methods when solving (\ref{eq:csdproblem}), we
provide below what we consider to be strong evidence for the
hypothesis that it is due to bad absolute condition numbers.  At the
very least we provide, in the form of a left preconditioner, an
effective way to significantly accelerate the numerical solution of
experimental design problems.

\section{The model problem and its conditioning}

We are led to our model problem by removing elements of the full
problem (\ref{eq:csdproblem}). The first simplification is the removal
of the external controls. As a consequence, the matrix $J$ is fixed,
and we have a problem in a form that is addressed by classical texts
(see e.g. \cite{Fedorov1972a}). Even this simplified problem is difficult
to analyze, as the assumptions imposed on $J$ in this theory are
rather weak. Informally speaking, once it has full rank, $J$ can be
any old matrix. Since this class of problems is very general, one can
expect to find counterexamples to any generalization of behavior
observed in practice. Our approach will thus be to restrict ourselves
to a model problem that is endowed with enough structure to
understand its badly-conditioned nature, and to devise a method of
curing it.

This problem is as follows. Given initial information in the form
of a variance-covariance matrix $\Sigma$ (which we, for simplicity, 
will consider to be a multiple $\alpha$ of the identity matrix) 
for the parameters of interest, our task is to choose between two
additional observations using a relaxed formulation. The idea of this
problem is to model an advanced stage of our experimental
design, in which, for the sake of argument, all weights except for
two are known to be either one or zero. The role of the parameter
$\alpha$ is to model the quality of the initial design that is the starting
point of our model problem. The smaller the $\alpha$, the better the
design that is to be improved.

The optimization problem we will study is thus
\begin{subequations}\label{eq:modelprob}
\begin{align}
\min_{w_{1},w_{2}}\mathrm{Tr}\left(\left[\alpha^{-1}I+w_1v_1v_1^T+w_2v_2v_2^T\right]^{-1}\right)
\intertext{subject to}
w_{1}+w_{2}=1\quad\text{and}\quad
0\leq w_{i}\leq 1, i=1,2.
\end{align}
\end{subequations}

\begin{theorem}\label{thm:ill}
Suppose that $v_1^Tv_2=0$ and $\|v_{1}\|=\|v_2\|=1$. Then the
solution of problem (\ref{eq:modelprob}) is $w_1=w_2=1/2$, and the
absolute condition number of the solution is
\begin{equation}
\label{eq:8}
\kappa_{\mathrm{abs}}=\frac{\left(1+\frac12\alpha\right)^3}{2\alpha^3}.
\end{equation}
\end{theorem}

The use of the absolute condition number is due to the fact that we
are searching for the zero of a derivative. The error in the ``data''
is thus a perturbation of the zero, which is only meaningful in an
absolute sense.

\begin{proof} %
  We write first $w_2=1-w_1$ to eliminate the equality constraint,
  and thus obtain a problem in one variable $w$. The function we want to
  minimize is then
\begin{equation}
\label{eq:1}
f(w):=\mathrm{Tr}\left(\left[\alpha^{-1}I+wv_1v_1^T+(1-w)v_2v_2^T\right]^{-1}\right).
\end{equation}
The orthogonality of $v_1$ and $v_2$ simplifies the application of the
Sherman-Morrison formula to obtain
\begin{equation}
\label{eq:5}
f(w):=2\alpha-\frac{w\alpha^2}{1+w\alpha}-\frac{(1-w)\alpha^2}{1+(1-w)\alpha}.
\end{equation}
Now we look for $w^{*}$ such that $f'(w^{*})=0$. Straight-forward
algebraic arguments yield $w^{*}=1/2$.

To investigate the effect of perturbations, we choose $\epsilon>0$,
and define the solution mapping
$g:(-\epsilon,\epsilon)\rightarrow\setR$ by
$f'(g(\epsilon))=\epsilon$. Of course, $g(0)=w^{*}$, and the condition
number of the problem $f'(w)=0$ is given by \cite{Demmel1987}
\begin{equation*}
\kappa_{\mathrm{abs}}=\frac{|g'(0)|}{|g(0)|}.
\end{equation*}
Using implicit differentiation on $f'(g(\epsilon))-\epsilon=0$ we
obtain the expression
\begin{equation}
\label{eq:7}
\kappa_{\mathrm{abs}}=\frac{1}{|f''(w^{*})||w^{*}|}.
\end{equation}

To finish the proof, we only need to compute $f''(w)$ and verify the
expression.
\end{proof}

Theorem \ref{thm:ill} suggests that as the optimization progresses and
our experimental design becomes better and better, that is, $\alpha$
becomes smaller, the absolute condition number of the problem will
increase. It will behave as $\alpha^{-3}$ for small
$\alpha$. Informally speaking, the bottom of the valley in which the
solution lies will become very flat, making it harder and harder
to choose between two additional rows of roughly the same size. The
model problem thus predicts the stagnation in the solution process,
which is precisely what often occurs in practice.

At least for the model problem, there is a surprisingly simple remedy
in the form of a {\em left preconditioner}.  A left preconditioner for
a given problem is a diffeomorphism in the dependent variables that
preserves the solution, and at the same time it improves the condition
number. Its name comes from the fact that it appears to the left of
the function that defines the unpreconditioned problem. The definition
we use here is an extension of the definition that is commonly used in
linear algebra (see e.g \cite{Saad2003}) to a nonlinear setting. The
next theorem suggests a possible choice of a left preconditioner for
the model problem.

\begin{theorem}\label{thm:uniform}
The problem
\begin{subequations}\label{eq:modmodelprob}
\begin{align}
\min_{w_{1},w_{2}}-\left\{\mathrm{Tr}\left(\left[\alpha^{-1}I+w_1v_1v_1^T+w_2v_2v_2^T\right]^{-1}\right)\right\}^{-2}
\intertext{subject to}
w_{1}+w_{2}=1\quad\text{and}\quad
0\leq w_{i}\leq 1, i=1,2.
\end{align}
\end{subequations}
has the same minimum as (\ref{eq:modelprob}).  For every $\alpha>0$
the absolute condition number of the minimum is $\kappa_{\mathrm{abs}}=2$.
\end{theorem}

\begin{proof} Repeat, with $\tilde{f}:=-f^{-2}$, the calculation of
  the condition number of Theorem \ref{thm:ill} 
\end{proof}

In light of the above, a possible way to achieve better convergence
when solving (\ref{eq:csdproblem}) is to take inspiration in Theorem
\ref{thm:uniform} and choose the left preconditioner
\[ h: (0,\infty) \to (- \infty, 0) \qquad h(z):=-z^{-2}. \]
As a consequence, we propose (and recommend) to solve the following
problem instead of \eqref{eq:csdproblem}:

\begin{subequations}\label{eq:csdproblempre}
\begin{align}
\min_{w,q}-\left\{\mathrm{Tr}\left(\left[J^{T}(q)W(w)J(q)\right]^{{-1}}\right)\right\}^{-2}\\
\intertext{subject to the constraints}
w\in\Omega(m_{\max},m)\quad\text{and}\quad q\in\Theta.
\end{align}
\end{subequations}
As we will see in the next section, this modification has the desired
effect of accelerating the convergence of SQP methods when solving
experimental design problems.

\section{Numerical experiments}
\label{sec:naexps}

In what follows, we will compare experimentally how well the SQP
solver behaves when solving the original problem (\ref{eq:csdproblem})
versus the preconditioned problem (\ref{eq:csdproblempre}). In the
first experiment, we tackle the task of optimizing a design when given
a prior information, but without external controls $q$. In the second
experiment, we keep $m_{\max}$ constant and vary the number of
candidate measurements to observe how the reformulation affects
performance as the size of the problem changes. The third numerical
experiment is a full nonlinear experimental design problem with
external controls on a model defined through a system of ordinary
differential equations. With it, we intend to illustrate the effect of
the reformulation on practical problems.

To make the results representative and reproducible, we programmed the
SQP solver as it is described in \cite{Nocedal2006}, with an augmented
Lagrangian penalty function as described in
\cite{Schittkowski82a}. This results in a reasonably robust and
effective solver.  The quadratic problems are solved using QPOPT
\cite{Gill1995corr}, and the Hessian is approximated using damped BFGS
updates \cite{Nocedal2006}. To avoid scaling issues as much as
possible, we use as an initial Hessian approximation a diagonal
matrix with the absolute values of the diagonal of the exact
Hessian. This retains the diagonal of the exact Hessian in the
unpreconditioned (convex) case, and assures positive definiteness 
in the preconditioned case.

We developed a Radau IIa solver of order 5 based on the ideas in
\cite{Hairer1999} that uses internal numerical differentiation
\cite{Bock1981} to compute accurate sensitivities of the differential
equation (here: 1st to 3rd order). All derivatives of algebraic
functions, including the inversion of the information matrix, were
computed using automatic differentiation
\cite{Griewank2008corrected,Griewank1989}, with a Common Lisp AD package
\cite{Fateman06buildingalgebra} extended for higher derivatives. The
inversion of the information matrix was done by directly computing
$J^TW(w)J$ and applying a Cholesky decomposition. To cope with
stability issues, we used quad-double arithmetic \cite{Hida2001corrected}. This
also ensures we can use the full range $[0,1]$ for the weights
\cite{Mommer2011}.

The problems without external controls are linear, and thus are
defined through the matrix $J$, which we generate randomly. We do this
by filling a matrix with random entries uniformly distributed in
$[-1,1]$, performing a singular value decomposition, and substituting
its singular values with exponentially decaying ones chosen to obtain
a condition number of $10^4$.

\subsection{Effect of prior information}

In this first experiment, we want to observe the behavior of SQP
methods on the sampling design problem for choosing 20 out of 50
possible measurement points:
\begin{equation}\label{eq:probwpreinfo}
\min_{w\in\Omega(20,50)}\Tr\left(\left[\alpha^{-1} I+ J^{T}W(w)J\right]^{-1}\right)
\end{equation}
with and without using the preconditioner. The matrix
$J\in\setR^{50\times7}$ is chosen at random as described before, but
we additionally normalize all its rows in the euclidean norm. We thus
obtain a nontrivial problem that is similar to our model, and where we
can observe the effect of prior information $\alpha I$ directly.

\begin{figure}
\begin{center}
\includegraphics[width=6.0cm]{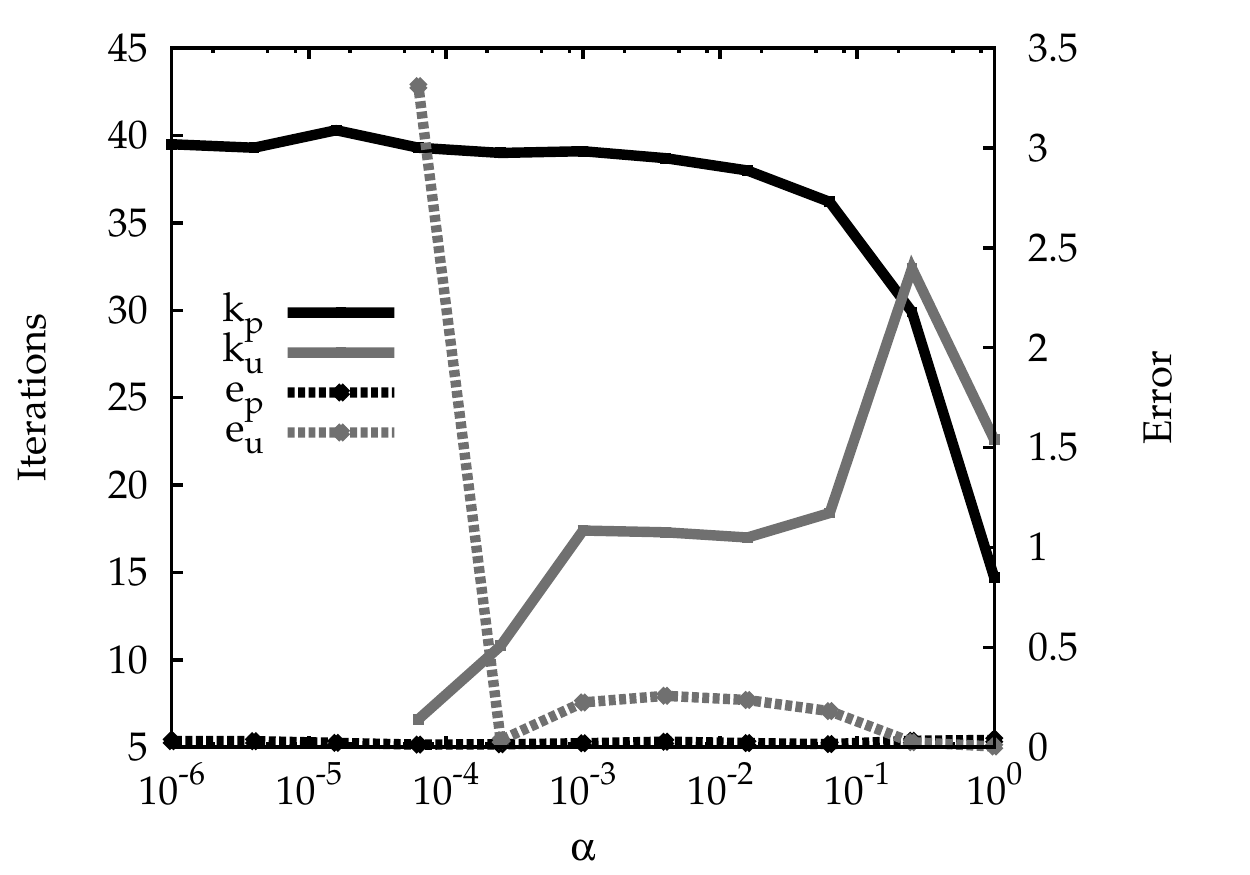}
\includegraphics[width=6.0cm]{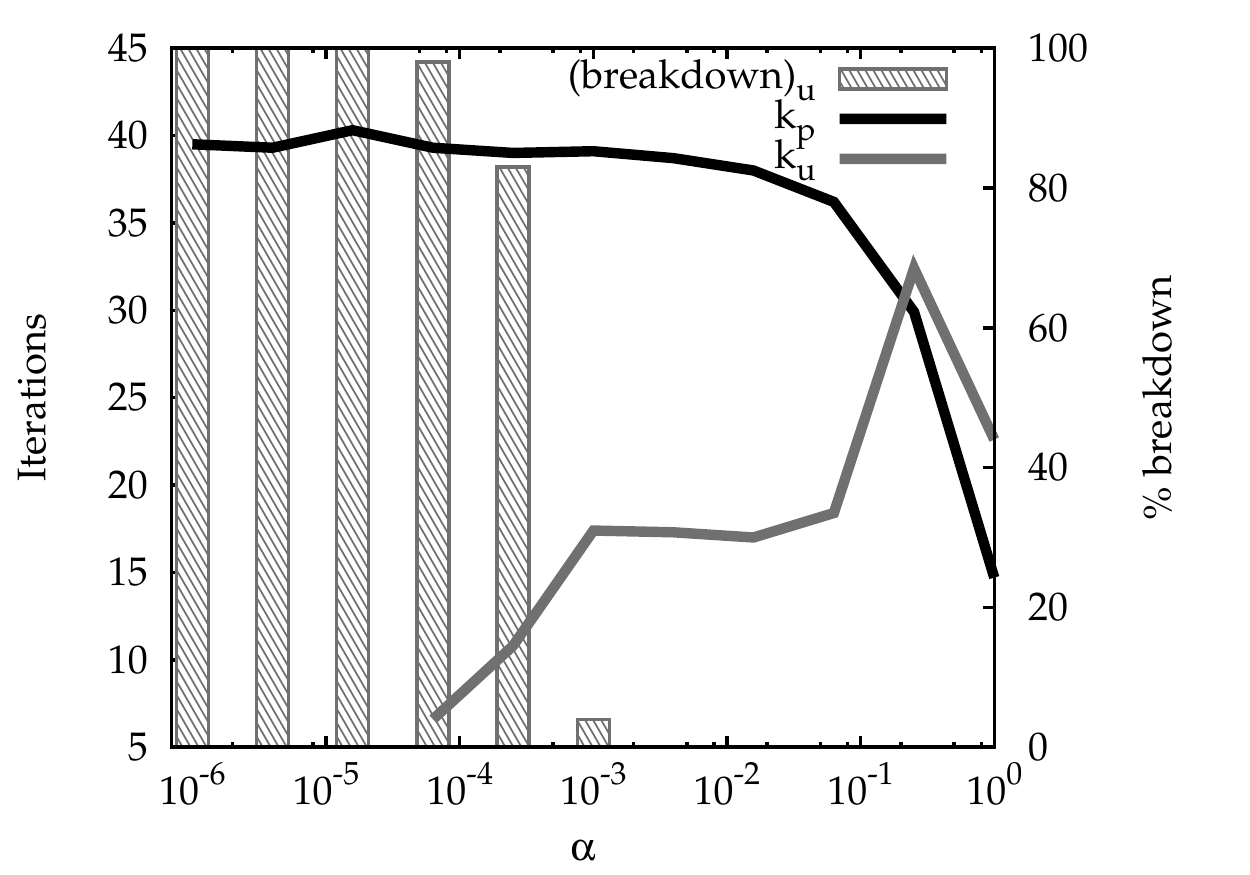}
\end{center}
\caption{Effect of prior information on stability and average
  iteration counts of the unpreconditioned and preconditioned problems
  (subscripts $u$ and $p$, respectively). See text for details.}
\label{fig:collapse}
\end{figure}

In Figure \ref{fig:collapse} we summarize the results of solving
problem (\ref{eq:probwpreinfo}) with and without preconditioning for
200 different matrices and 11 different values of $\alpha$, chosen 
equidistantly on a logarithmic scale in $[10^{-6},1]$. 
As $\alpha$ decreases, and thus the prior
information becomes better and better, we see that the iteration count
of the preconditioned variant stabilizes to around $40$, while the
sensitivity to the initial guess is about the same for each
$\alpha$. This distance was estimated by comparing the solution we
obtained using two different starting guesses for the iteration. One
that has the first $20$ components equal to $1$, and the rest is zero,
and one starting guess that is the reverse. Since the problem is
convex, we can use this sensitivity as a proxy for the distance to the
exact solution. This value tends to be a lot larger for the
unpreconditioned problem, which also becomes very hard to solve for
small $\alpha$. On the right of Figure \ref{fig:collapse} we have
plotted the percentage of problems that could not be solved because
the quadratic solver reached its default iteration limit, something
that never occurred with the preconditioned formulation.

\subsection{Effect of problem size}

\begin{figure}
\includegraphics[width=12cm]{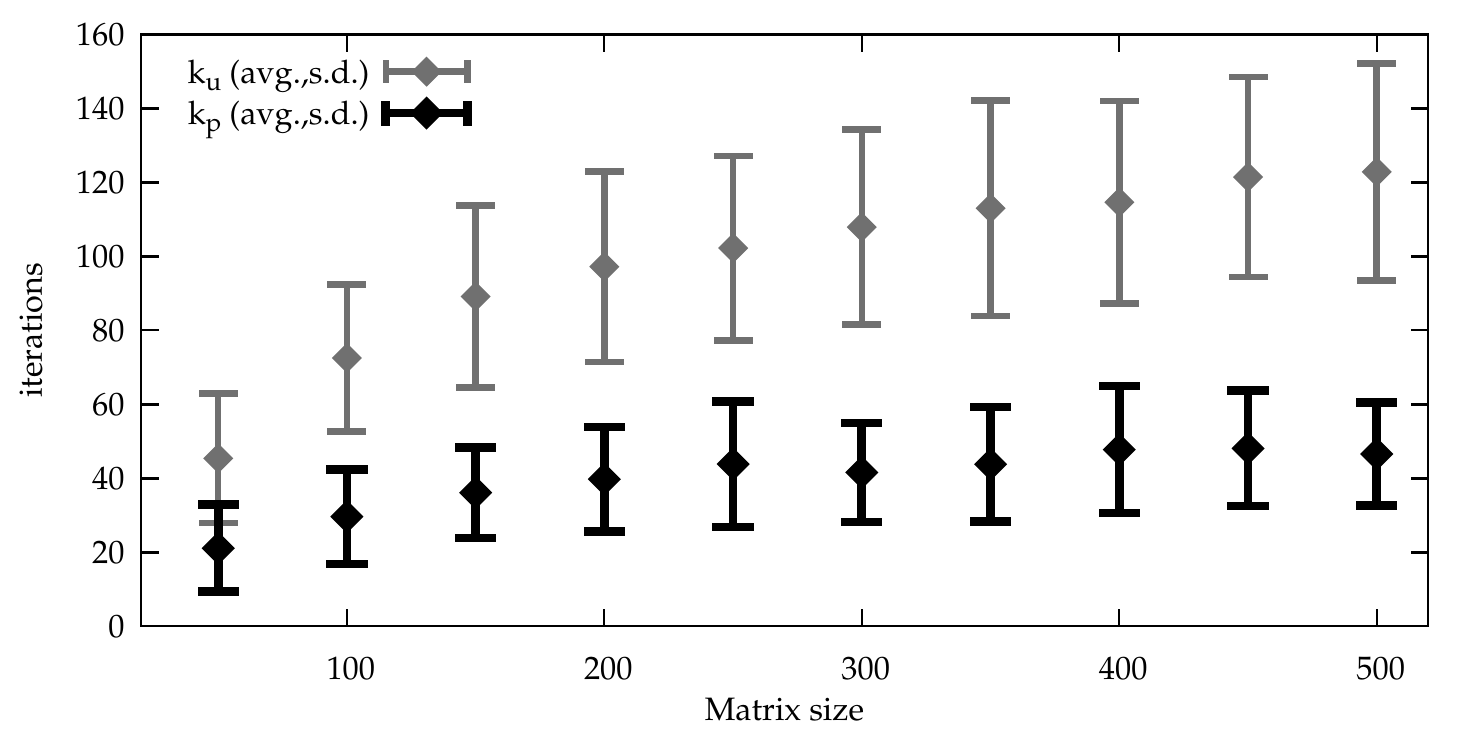}
\caption{Average iterations of the SQP method solving the
  unpreconditioned and preconditioned problem (subscripts $u$ and $p$,
  respectively).}
\label{fig:avgits}
\end{figure}

Now we consider matrices $J$ of size $50n\times7$ for
$n=1,2,\ldots,10$. For each size we generated $200$ matrices, and
solved the problem
\begin{equation}
\label{eq:samplingdesprob}
\min_{w\in\Omega(20,50n)}\Tr\left(\left[J^{T}W(w)J\right]^{-1}\right)
\end{equation}
with and without preconditioner. We plot average iteration
counts in Figure \ref{fig:avgits}, with error bars giving the standard
deviations. Again we observe a significant improvement, on average, of
the iteration counts for each problem. We also observe that,
for the unpreconditionend problem, the iteration count increases with
the problem size $n$, which is
likely due to the fact that the value of the minimum becomes smaller
with matrix size, so that the effect predicted by Theorem
\ref{thm:ill} increases. Whether the iteration counts increase or not
for the preconditioned variant is not possible to tell conclusively
from this experiment, but we conjecture that it does.

\subsection{Full nonlinear problem}

Our original motivation was to find an advantageous formulation for
full nonlinear experimental design as it is relevant for applications in
engineering. Thus we consider in our next experiment an experimental
design problem defined on a system of nonlinear differential
equations. For simplicity, we choose the FitzHugh-Nagumo
\cite{FitzHugh1961,Nagumo1962} model, which is given by the system
\begin{align}
  \dot x_1 &= x_1 - z x_1^3 - x_2 + I   \qquad~ 
  &  x_1(t_0) = x_{0,1}, \\
  \dot x_2 &= a \cdot (x_1 + b + cx_2) \qquad~ 
  &  x_2(t_0) = x_{0,2}. 
\end{align}
The task is to find a experimental design to estimate the four fixed
parameters $z=0.25$, $a=0.02$, $b=0.7$, and $c=-0.8$ using the values
of $I$, $x_{0,1}$, and $x_{0,2}$ as controls, which are expected to
satisfy the constraints
\begin{equation}
\label{eq:2}
-5<x_{0,1}<5,\quad-5<x_{0,2}<5,\quad -1<I<0.5.
\end{equation}
At most $30$ measurements should be taken. A measurement is the
reading of the value of any of the two variables at any of the times
$t_i=5i$, $i=1,2,\ldots,100$.

\begin{table}
\center
\begin{tabular}{|l|l|l||l|l|l|l|}
\hline
  score & $\langle k_{p} \rangle$  & $\langle k_{u} \rangle$ &
  $\langle k_{u}/k_{p}\rangle$ & $\sigma$  \\
\hline
  5:0 & 46.0 & 260.8 & 3.9 & 1.7 \\
\hline
\end{tabular}
\caption{Performance statistics for the FitzHugh-Nagumo example. The
  values of $k$ denote iteration counts, the subscripts $u$ and $p$
  indicate that they refer to unpreconditioned and preconditioned
  variants. The angled brackets indicate average, and $\sigma$ stands
  for the standard deviation of the speed-up factor $\langle k_{u}/k_{p}\rangle$ }
\label{tab:fhn}
\end{table}

The initial guess for the weights is that all of them are equal. We
choose the initial guess for the optimal controls randomly to be able
to assess typical behavior as much as possible. Since it is possible
to choose the controls in such a way that the least-squares problem is
singular, we only use those for which
$\Tr([J^{T}(q_{0})J(q_{0})]^{-1})\leq 100$. The optimal values are
usually three to four orders of magnitude lower than this upper limit.

\begin{figure}
\includegraphics[width=12cm]{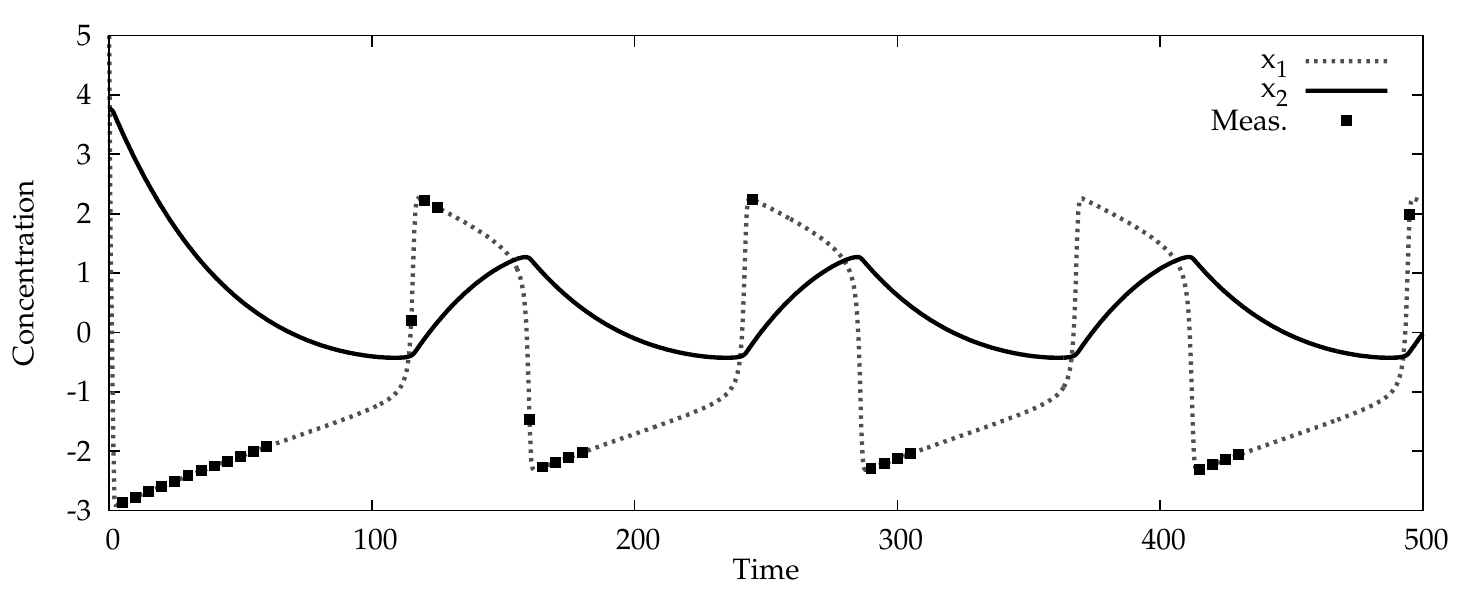}
\caption{Optimal design for the FitzHugh-Nagumo system}
\label{fig:fhnplan}
\end{figure}

We repeated the experiment five times, and summarize the results in
table \ref{tab:fhn}. We note again that the preconditioned formulation
leads to important gains in convergence speed. Here we have that the
problem is not convex for either formulation, and one observes that
there exist many local minima. For a given starting value, the two
variants may or may not converge to the same one. In Figure
\ref{fig:gure} we can observe the convergence behavior of a typical
run in terms of the length of the search direction. For completeness,
we include in Figure \ref{fig:fhnplan} a solution of this optimization
problem ($I=0.37$, $x_{0,1}=5$, $x_{0,2}=3.76$), together with the
optimal measurement points.

\section{Final remarks and outlook}

In this article, we have studied the convergence problems of SQP
methods when solving experimental design problems. Our findings 
suggest that bad absolute condition numbers are indeed the cause of
slow convergence.

We work around the generality of the problem by studying a carefully
chosen model problem involving the $A$-criterion. From the
understanding gained about the conditioning of this model problem we
derive a transformation that guarantees constant absolute condition
numbers in this particular case.

We then provide strong experimental evidence that this transformation
yields significant improvements in the convergence behavior of SQP
methods when applied to relevant settings, where we observe that the
improvement is more significant for larger problems. The
transformation itself does not introduce any noticeable additional
computational cost.

\begin{acknowledgements}
  This work was supported by the German Ministry of Education and
  Research (BMBF) under grant ID: 03MS649A, and the Helmholtz
  association through the SBCancer program.
\end{acknowledgements}

\bibliographystyle{spmpsci}
\bibliography{artbib,agbock}

\end{document}